\documentclass[centertags,leqno]{article}

\usepackage{amsmath}
\usepackage{amssymb}
\usepackage{latexsym}
\usepackage{amsxtra}
\usepackage{amscd}
\usepackage{theorem}
\usepackage{epic}
\usepackage{eepic}

\theoremstyle{plain}

{\theorembodyfont{\slshape}

        \newtheorem{thm}{Theorem}
        
        \newtheorem{lem}[thm]{Lemma}
        \newtheorem{prop}[thm]{Proposition}
}

{\theorembodyfont{\rmfamily}

        \newtheorem{prob}[thm]{Problem}
        \newtheorem{rem}[thm]{Remark}

}

\renewcommand{\em}{\sl}

\newcommand{\proof}{{\bf Proof:\ }}
\newcommand{\Endproof}{\hspace*{\fill} $\Box$ \vspace{1ex} \noindent }

\makeatletter
\renewcommand{\subsection}{\@startsection{subsection}{2}%
        {\z@}{-3.25ex plus -1ex minus-.2ex}{-1em}{\bf}}
\makeatother

\newcommand{\PP}{\mathbb{P}}
\newcommand{\ZZ}{\mathbb{Z}}

\newcommand{\QQ}{\mathbb{Q}}

\newcommand{\FF}{\mathbb{F}}
\renewcommand{\AA}{\mathbb{A}}

\newcommand{\M}{\mathcal{M}}
\newcommand{\HH}{\mathcal{H}}

\newcommand{\X}{\mathcal{X}}
\newcommand{\Y}{\mathcal{Y}}
\newcommand{\Z}{\mathcal{Z}}
\newcommand{\PSL}{\mathop{\rm PSL}}

\newcommand{\Spec}{{\rm Spec\;}}
\newcommand{\Ind}{{\rm Ind}}

\newcommand{\alg}{^{\rm alg}}
\newcommand{\stab}{^{\rm stab}}

\newcommand{\gen}[1]{\mathopen\langle#1\mathclose\rangle}

\newcommand{\zz}{\text{$2$-$2$}}

\title{Alternating groups as monodromy groups in positive characteristic}

\author{Stefan Wewers and Irene I.\ Bouw}
\date{}

\begin{document}

\maketitle
\begin{abstract}
Let $X$ be a generic curve of genus $g$ defined over an algebraically
closed field $k$ of characteristic $p\geq 0$. We show that for $n$
sufficiently large there exists a tame rational map $f:X\to \PP^1_k$
of degree $n$ with monodromy group $A_n$. This generalizes a result of
Magaard--V\"olklein to positive characteristic.
\end{abstract}

\section{Introduction} \label{intro}

\subsection{} \label{intro1}

Let $k$ be an algebraically closed field of characteristic
$p\geq 0$. We denote by $k_0$ the prime field of $k$ (i.e.\ 
$k_0=\QQ$ if $p=0$, $k=\FF_p$ if $p>0$). For $g\geq 0$, we
write $\M_g$ for the (coarse) moduli space of curves of
genus $g$ in characteristic $p$. This is a smooth,
quasi-projective and geometrically irreducible variety over
$k_0$, see \cite{DeligneMumford}. Its dimension is $3g-3$
(resp.\ $g$) for $g\geq2$ (resp.\ for $g=0,1$).

Let $X$ be a smooth projective curve of genus $g$ over $k$.
It corresponds to a $k$-rational point $x:\Spec k\to\M_g$.
We say that $X$ is {\em generic} if the image of $x$ is
Zariski dense in $\M_g$. 

Let $f:X\to\PP^1_k$ be a non-constant rational function on $X$. We say
that $f$ is {\em tame} if the extension of function fields $k(X)/k(f)$
is separable and at most tamely ramified. The degree of the extension
$k(X)/k(f)$ is called the {\em degree} of $f$.  The Galois group of
the Galois closure of $k(X)/k(f)$ is called the {\em monodromy group}
of $f$. The aim of this note is to prove the following theorem.

\begin{thm} \label{thm0}
  Let $g\geq 0$ and $n\geq 3$. Let $X$ be a generic curve of
  genus $g$, defined over an algebraically closed field $k$
  of characteristic $p\geq 0$. Then the curve $X$ admits a
  tame rational function $f:X\to\PP^1_k$ of degree $n$ with
  alternating monodromy group, in each of the following
  cases:
  \begin{enumerate}
  \item[(a)]
    If $p\neq 2,3$ and $n\geq \max\,(g+3,2g+1)$.
  \item[(b)] If $p=2$, $n\geq \max\,(g+3,2g+1)$ and
    $n+g$ is odd.
  \item[(c)] If $p=2$, $n\geq \max\,(g+6,2g+3)$ and $n+g$ is even.
  \item[(d)] If $p=3$ and $n\geq\max\,(7,g+6,2g+1)$.
  \end{enumerate}
\end{thm}

In characteristic $p=0$, this has recently been proved by Magaard and
V\"olklein  \cite{VoelkleinMagaard03} (except for the case
$(g,n)=(2,5)$). The cases $g=1$ and $p=0$ had been proved earlier by
Fried, Klassen and Kopeliovich \cite{FriedKlaKop00}. Also in
characteristic $0$, Artebani and Pirola \cite{ArtebaniPirola} have
shown that {\em every} curve admits a tame rational function of degree
$n$ with alternating monodromy group, provided that $n\geq 12g+4$. In
characteristic $p\not=3$, Schr\"oer \cite{Schroer02} has proved that
for every $g\geq 0$ there exists {\em some} curve of genus $g$
admitting a rational function with monodromy group $A_n$, for certain
values of $n$.

Note that there are classical analogs of these results
for the symmetric group. In fact, one knows that every curve of genus
$g\geq 2$ over an algebraically closed field of characteristic
$p\not=2$ admits a tame rational function of degree $n$ with monodromy
group $S_n$, provided that $2n-2\geq g$ (\cite{Fulton69}). However,
in characteristic $2$, it seems to be unknown whether every
curve admits a tame rational function at all -- no matter what the
degree and the monodromy group are.

Our results are somewhat more precise than Theorem \ref{thm0}
above: the tame cover $f:X\to\PP^1_k$ can be chosen in such a way that
its inertia groups are all generated by $3$-cycles (except if $p=3$ or
if $p=2$ and $n+g$ is odd). If $p\not=2,5$ and $n\geq \max\,(g+6,2g+1)$,
the cover $f$ can be chosen in such a way that its inertia groups are
all generated by double transpositions (with the possible exception
$(p,g,n)=(3,0,6)$). See the statement of Theorem \ref{thm1}.

For $p\not=2,3$ our results are optimal, in the following
sense. If $n<g+3$ or $n<2g+1$ then the generic curve of
genus $g$ does {\em not} admit a tame rational function of degree
$n$ with alternating monodromy group. (For $g\geq 3$, this is
shown in \cite{VoelkleinMagaard03}, see the proof of Theorem
3.3 in {\em loc.cit}. For $g<3$, one can use similar arguments.)

For $p=2$ or $p=3$, it is unlikely that our results are optimal. For
instance, we believe (but were not able to show) that the generic
curve of genus $1$ in characteristic $2$ admits a tame rational
function of degree $5$ with alternating monodromy. If this were the
case, then the situation in characteristic $2$ would be the same as in
characteristic $0$, i.e.\ the condition `$n+g$ odd' in Theorem
\ref{thm0} (b) would be unnecessary (except for the case
$(g,n)=(0,4)$, which does {\em not} occur in characteristic $2$). This
would then also be an optimal result.  But at the moment, there are
infinitely many pairs $(g,n)$ which, to our knowledge, may or may not
occur in characteristic $2$. In characteristic $3$, there are only
finitely many such cases. See Section \ref{OpenProblems} for a list of
all open cases.

\subsection{} \label{intro3}

In order to show the existence of rational functions with alternating
monodromy, the authors of \cite{VoelkleinMagaard03} use {\em Hurwitz
spaces}, i.e.\ moduli space for covers of $\PP^1$. In fact, the main
result of \cite{VoelkleinMagaard03} is stated in the following form:
the natural map $\HH_{r,n}\to\M_g$ from a certain Hurwitz space to the
moduli space of curves (which maps the isomorphism class of a cover
$f:X\to\PP^1$ to the isomorphism class of the curve $X$) has a dense
image.  In \cite{FriedVoelklein91}, Hurwitz spaces are constructed in
characteristic $0$, using tools from topology and the theory of
complex analytic functions. From this point of view, it seems
difficult to extend the results of \cite{VoelkleinMagaard03} to
positive characteristic.  However, it is shown in \cite{Fulton69} and
\cite{diss} that Hurwitz spaces can also be constructed in a purely
algebraic way and therefore make sense in positive characteristic,
too.  Moreover, for almost all primes $p$ a given Hurwitz space has
good reduction from characteristic $0$ to characteristic $p$ (in the
case of the Hurwitz space $\HH_{r,n}$ used in
\cite{VoelkleinMagaard03}, this is true for all primes $p>n$).  Using
this good reduction result, it is then easy to extend the results of
\cite{VoelkleinMagaard03} to characteristic $p$, provided that $p>n$.

However, for small primes $p$ this kind of argument does not work. For
instance, if $p\leq n$ the Hurwitz space $\HH_{n,r}$ may have bad
reduction to characteristic $p$.  This makes it very difficult to
study $\HH_{n,r}\otimes\FF_p$, using results on
$\HH_{n,r}\otimes\QQ$. The same phenomenon prevented Fulton
\cite{Fulton69} from proving the irreducibility of $\M_g\otimes\FF_p$
for $p\leq 2g+1$.

Luckily, it turns out that most of the arguments in
\cite{VoelkleinMagaard03} (for instance, the induction step
from $(g,n)$ to $(g+1,n+1)$) can be carried out in a purely
algebraic way. In fact, Hurwitz spaces are not really
essential; the key results that are used are the geometry of
the boundary of $\M_g$ and the deformation theory of tame
covers, developed in \cite{SGA1}. Both these tools are
algebraic and work in any characteristic.  The situation is
different for the proof of Theorem \ref{thm0} for small
values of $g$ and $n$, where the induction process starts.
Here the arguments in \cite{VoelkleinMagaard03} are based on
Riemann's Existence Theorem and do not carry over to
characteristic $p$ if $p$ is small.  This is the reason why
the statement of Theorem \ref{thm0} is weaker for $p=2,3$.

\vspace{3ex} The authors would like to thank Gerhard Frey
for stimulating this work and Helmut V\"olklein and the
referee for useful comments on an earlier version of this manuscript.

\section{The two induction steps} \label{ind}

\subsection{Deformation of admissible covers} \label{deform}

Let $k$ be an algebraically closed field. We choose a compatible
system $(\zeta_n)$ of $n$th roots of unity in $k$, where $n$ runs over
all natural numbers prime to the characteristic of $k$. Let $t$ denote
a transcendental element over $k$. Our goal is to construct tame
covers $f_\eta:X_\eta\to\PP^1_{k((t))}$ over the field $k((t))$ by
deforming a given tame cover $f_s:X_s\to Z_s$ between singular
curves. The standard general reference for the deformation theory of
tame covers is of course \cite{SGA1}. For the particular results that
we use, we refer to \cite{Saidi97}, \cite{HarbaterStevenson} or
\cite{admissible}.

Let $G$ be a finite group.  For $i=1,2$, we have a subgroup
$G_i\subset G$, a tame $G_i$-Galois cover $h_i:Y_i\to\PP^1_k$ between
smooth projective curves over $k$, and a closed point $y_i\in Y_i$. We
assume that the datum $(G,G_i,h_i,y_i)$ has the following properties.
First, we assume that $G$ is generated by its subgroups $G_1$ and
$G_2$. For the second condition, let $g_i\in G_i$ be the canonical
generator of the stabilizer of $y_i$, with respect to $(\zeta_n)$. (An
element $g\in G_i$ with $g(y_i)=y_i$ is called a {\sl canonical generator} if
there exists a formal parameter $u$ at $y_i$ such that
$g^*u=\zeta_n\cdot u$, with $n={\rm ord}(g)$.) Then we demand that
$g_1=g_2^{-1}$. We denote by $n_0$ the order of $g_1=g_2^{-1}$.

Given $(G,G_i,h_i,y_i)$ satisfying the above conditions, we construct
a tame $G$-Galois cover $h_s:Y_s\to Z_s$ between semistable curves
over $k$, as follows. We set
\[
   Y_s \;:=\; \Big(\,\Ind_{G_1}^G(Y_1)\;{\textstyle \coprod}\;
           \Ind_{G_2}^G(Y_2)\,\Big)/_\sim,
\]
where $\sim$ denotes the following equivalence relation. A point $y\in
\Ind_{G_1}^G(Y_1)$ is equivalent to a point $y'\in \Ind_{G_2}^G(Y_2)$
if and only there exists an element $g\in G$ with $y=g(y_1)$ and
$y'=g(y_2)$. It is easy to see that the set $Y_s$ is naturally
equipped with the structure of a semistable curve over $k$ and with a
$k$-linear action of $G$. The curves $Y_1\subset Y_s$ and $Y_2\subset
Y_s$ are irreducible components of $Y_s$, with stabilizer $G_1$ and
$G_2$, respectively. Moreover, the points $y_1\in Y_1$ and $y_2\in
Y_2$ correspond to the same (singular) point of $Y_s$. We define
$Z_s:=Y_s/G$ as the quotient scheme. The scheme $Z_s$ is a semistable
curve over $k$, with two irreducible components $Z_1,Z_2$ which meet
in one points. The components $Z_i$ can be identified with $\PP^1_k$,
via the covers $h_i$. The natural map $h_s:Y_s\to Z_s$ is a {\em tame
admissible cover} (\cite{admissible}).

It is easy to construct a semistable curve $\Z$ over $\Spec k[[t]]$
with special fiber $Z_s$ and with generic fiber
$Z_\eta=\PP^1_{k((t))}$, satisfying the following additional property:
the complete local ring of $\Z$ at the singular point of the special
fiber is isomorphic to $k[[t,u,v\mid uv=t^{n_0}]]$. Let
$\bar{\tau}_1,\ldots,\bar{\tau}_s\in Z_1=\PP^1_k$ denote the branch
points of $h_1$ distinct from $h_1(y_1)$, and
$\bar{\tau}_{s+1},\ldots,\bar{\tau}_r$ the branch points of $h_2$
distinct from $h_2(y_2)$. Lift these points to $k[[t]]$-rational
points $\tau_1,\ldots\tau_r$ of $\Z$. By \cite{admissible}, Theorem
3.1.1, there exists a tame $G$-Galois cover $h:\Y\to\Z$ between
semistable curves over $k[[t]]$, \'etale over
$\Z-\{\tau_1,\ldots,\tau_r\}$ whose special fiber is equal to the
cover $h_s$. The generic fiber $h_\eta:Y_\eta\to\PP^1_{k((t))}$ is a
tame $G$-Galois cover between smooth projective curves, with $r$
branch points. We say that $h_\eta$ is a {\em smooth Galois cover}
associated to $(G,G_i,h_i,y_i)$.

Suppose that $G\subset S_n$ is a transitive permutation
group on $n$ letters, and let $H\subset G$ denote the
stabilizer of $1$ in $G$. Let $\X:=\Y/H$ denote the quotient
scheme. Then $\X$ is a semistable curve over $k[[t]]$, with
generic fiber $X_\eta=Y_\eta/H$ and special fiber
$X_s=Y_s/H$ (for the second equality we have used the fact
that the cover $h_s:Y_s\to Z_s$ is separable). Moreover,
$h_\eta$ is the Galois closure of the cover
$f_\eta:\X_\eta\to\PP^1_{k((t))}$. Therefore, the cover
$f_\eta$ is tame, with monodromy group $G$. We say that
$f_\eta$ is a {\em smooth cover} associated to $(G\subset
S_n,G_i,h_i,y_i)$.

It is easy to describe the special fiber $X_s$ and the admissible
cover $f_s:X_s\to Z_s$ in terms of the datum $(G\subset
S_n,G_i,h_i,y_i)$. For instance, the irreducible components of $X_s$
lying above the component $Z_i\subset Z_s$ are in bijection with the
orbits of the  $G_i$-action on $\{1,\ldots,n\}$. For each orbit
$O\subset\{1,\ldots,n\}$, the restriction of $f_s$ to the component
$X_O\subset X_s$ is isomorphic to the quotient cover of $h_i:Y_i\to
Z_i=\PP^1_k$ corresponding to the stabilizer of some element of $O$.

\subsection{A useful lemma} \label{2pt}

Let $X$ be a generic curve of genus $g\geq 1$, defined over an
algebraically closed field $k$ of characteristic $p\geq 0$. Let
$f:X\to\PP^1_k$ be a tame cover, with monodromy group $G$ and with
branch points $t_1,\ldots,t_r$, $r\geq 3$. Without loss of generality,
we may assume that $t_1=0$, $t_2=1$ and $t_3=\infty$. We say that $f$
has {\em generic branch points} if the branch points $t_4,\ldots,t_r$,
considered as elements of $k$, are algebraically independent over the
prime field $k_0\subset k$.

The following lemma will be useful in the proof of Theorem
\ref{thm0}. 

\begin{lem} \label{2ptlem}
  Let $k'/k$ be a field extension, $z:\Spec k'\to\PP^1_k$ a generic
  point and $x_1,x_2:\Spec k'\to X$ two distinct $k'$-rational points
  with $z=f(x_1)=f(x_2)$. Suppose that one of the following conditions
  holds.
  \begin{enumerate}
  \item[(a)]
    We have $g=1$ and the monodromy group $G$ is doubly transitive.
  \item[(b)]
    We have $g\geq 2$ and $f$ has $r>3g$ generic branch points. 
  \end{enumerate}
  Then $(X\otimes k';x_1,x_2)$ is a generic two-pointed curve, i.e.\
  its classifying map $\Spec k'\to\M_{g,2}$ has a Zariski dense image.
\end{lem}

\proof Assume that Condition (a) of the Lemma holds. Let
$\varphi_1:\Spec k'\to\M_{1,1}$ (resp.\ $\varphi_2:\Spec
k'\to\M_{1,2}$) denote the classifying map of the
one-pointed curve $(X\otimes k',x_1)$ (resp.\ of the
two-pointed curve $(X\otimes k';x_1,x_2)$).  Recall that
$\M_{1,1}\cong\AA^1$ and that $\varphi_1$ is simply the
$j$-invariant of the elliptic curve $(X\otimes k',x_1)$. It
is well known that the $j$-invariant of an elliptic curve
depends only on the underlying curve and not on the
distinguished point (i.e.\ $\M_{1,1}\cong\M_1$). Hence we
may regard $\varphi_1$ as an element of $k$. Since $X$ is
generic, $\varphi_1$ is transcendental over $k_0$. It
follows that the image of $\varphi_2$ is contained in the
fiber $\M_{1,2}\otimes k_0(\varphi_1)$ of the natural map
$\M_{1,2}\to\M_{1,1}$ over $\Spec k_0(\varphi_1)$. To prove
the lemma, it suffices to show that the image of $\varphi_2$
is Zariski dense in $\M_{1,2}\otimes k_0(\varphi_1)$.  Since
$\M_{1,2}\otimes k_0(\varphi_1)$ is $1$-dimensional, it even
suffices to find a $k$-rational place $v$ of $k'$ such that
the image of $v$ under $\varphi_2$ lies on the boundary of
$\M_{1,2}\otimes k_0(\varphi_1)$. In other words, we require
$x_1(v)=x_2(v)$ (here $x_i(v)$ denotes the $k$-rational
point on $X$ obtained by `specializing' $x_i$ at $v$).

Let $S\subset X\times_k X$ denote the locus of pairs $(x,y)$ with
$f(x)=f(y)$. It is a $1$-dimensional closed subset, and the natural
map $S\to\PP^1_k$ is finite. The pair $(x_1,x_2)$ is a $k'$-rational
point on $S$ which maps to the generic point of $\PP^1_k$. Therefore,
the closure of the image of $(x_1,x_2)$ is an irreducible component
$S'$ of $S$, distinct from the diagonal $\Delta\subset X\times_k X$.
We have to show that $S'$ has nonempty intersection with $\Delta$. The
assumption that $G$ is doubly transitive implies that $S-\Delta$ is
irreducible. Hence we have $S'=\overline{S-\Delta}$. On the other
hand, since the cover $f:X\to\PP^1_k$ is not \'etale, $\Delta$ has
nonempty intersection with $S'=\overline{S-\Delta}$. This proves Lemma
\ref{2ptlem}, assuming Condition (a).

Assume now that Condition (b) holds. Let $k_1\subset k$ denote the
minimal algebraically closed subfield of $k$ over which the curve $X$
can be defined. Since $X$ is generic, $k_1$ is isomorphic to the
algebraic closure of the function field of $\M_g$. In particular,
${\rm tr.deg}(k_1/k_0)=3g-3$. Write $X=X_1\otimes_{k_1}k$. 

By assumption, the branch points $t_4,\ldots,t_r$ of $f$, considered
as elements of $k$, are algebraically independent over
$k_0$. Moreover, we have $r-3>3g-3={\rm tr.deg}(k_1/k_0)$. It follows
that the cover $f:X_1\otimes k\to\PP^1_k$ is not isotrivial, with
respect to the extension $k/k_1$. More precisely, the subfield
$k(f)\subset k(X)$ cannot be generated over $k$ by an element $f_1\in
k_1(X_1)$. 

Let $f_V:X_1\times V\to\PP^1_V$ be a model of $f$ over $V$, where $V$
is a variety over $k_1$ (this means that $f$ is the pullback of $f_V$
via a generic point $\Spec k\to V$). For each closed point $v\in V$ we
obtain, by specializing $f_V$, a tame cover
$f_v:X_1\to\PP^1_{k_1}$. Since $f$ is not isotrivial, there exist
infinitely many points $v\in V$ which give rise to pairwise weakly
non-isomorphic covers $f_v$. Now \cite{VoelkleinMagaard03}, Lemma 2.3,
shows that the locus $S\subset X_1\times_{k_1} X_1$ of pairs of points
$(x,y)$ which satisfy $f_v(x)=f_v(y)$ for all $v\in V$ is Zariski
dense in $X_1\times_{k_1} X_1$. On the other hand, since
$S\to\PP^1_{k_1}$ is finite and $f(x_1)=f(x_2)=z$ is generic, the pair
$(x_1,x_2)$ is a generic point of $S$. This proves Lemma \ref{2ptlem},
assuming Condition (b).  \Endproof

\begin{rem} 
  Lemma \ref{2ptlem} is very similar to \cite{VoelkleinMagaard03},
  Lemma 2.4, at least in the case $g\geq 2$. In the case $g=1$, Lemma
  \ref{2ptlem} improves \cite{VoelkleinMagaard03}, Lemma 2.4, by
  removing the assumption that $f$ depends on more parameters than
  $X$.  This is used in the proof of Theorem \ref{thm0} for the case
  $(g,n)=(2,5)$.
\end{rem}

\subsection{The first induction step} \label{ind1}

As before, we let $X$ be a generic curve of genus $g\geq 0$, defined
over an algebraically closed field $k$ of characteristic $p\geq 0$. 
The following proposition will serve as an induction step in the proof of
Theorem \ref{thm0}, from $(g,n)$ to $(g,n+2)$. 

\begin{prop} \label{indprop1}
  Assume that there exists a tame cover $f:X\to\PP^1_k$ of degree
  $n\geq 3$, with monodromy group $A_n$ and with $r$ branch
  points. Let $k':=k((t))\alg$ and $X':=X\otimes_k k'$. Then there
  exists a tame cover $f':X'\to\PP^1_{k'}$ of degree $n+2$, with
  monodromy group $A_{n+2}$ and with $r+2$ branch points.
\end{prop}

\proof Let $h_1:Y_1\to\PP^1_k$ be the Galois closure of $f$.
By assumption, the Galois group of $h_1$ is $G_1:=A_n$,
which we view as a subgroup of $G:=A_{n+2}$, in the obvious
way. Let $y_1\in Y_1$ be any closed point where $h_1$ is
unramified. If $p\not=3$ (resp.\ if $p=3$), we let
$h_2:Y_2\to\PP^1_k$ be a cyclic cover of order $3$ (resp.\ 
of order $2$), ramified at two points. We identify the
Galois group of $h_2$ with the subgroup $G_2\subset G$
generated by the $3$-cycle $(n,n+1,n+2)$ (resp.\ by the
double transposition $(n-1,n+1)(n,n+2)$). Note that
$G=\gen{G_1,G_2}$ in both cases. We also choose a point
$y_2\in Y_2$ where $h_2$ is unramified. By the construction
of \S \ref{deform}, there exists a tame $G$-Galois cover
$h':Y'\to\PP^1_{k'}$ which lifts the datum
$(G,G_i,h_i,y_i)$. Let $H\subset G$ be the stabilizer of
$1$. Set $X':=Y'/H$ and let $f':X'\to\PP^1_{k'}$ denote the
natural map. By construction, $f'$ is a tame cover of degree
$n+2$, with monodromy group $A_{n+2}$ and with $r+2$ branch
points.

It remains to prove that $X'\cong X\otimes_k k'$. The Riemann--Hurwitz
formula shows that the genus of $X'$ is equal to $g$. Moreover, by
construction, one of the components of the stable reduction of $X'$ is
isomorphic to $X$. It follows that $X'$ has good reduction and
specializes to a generic curve of genus $g$. Therefore, $X'$ is itself
generic, hence $X'\cong X\otimes_k k'$.  \Endproof

\subsection{The second induction step} \label{ind2}

For the proof of Theorem \ref{thm0}, we need another induction step,
going from $(g,n)$ to $(g+1,n+1)$. As in the last subsection, $X$
is a generic curve of genus $g\geq 0$ over $k$.

\begin{prop} \label{indprop2}
  Assume that there exists a tame cover $f:X\to\PP^1_k$ of degree
  $n\geq 3$, with monodromy group $A_n$ and with $r$ branch
  points. Assume, moreover, that one of the following conditions hold:
  \begin{enumerate}
  \item
    $g\leq 1$, 
  \item
    $g\geq 2$ and $r>3g$. 
  \end{enumerate}
  Then there exists a field extension $k'/k$, a generic curve $X'$ of
  genus $g+1$ over $k'$ and a tame cover $f':X'\to\PP^1_{k'}$ of
  degree $n+1$, with monodromy group $A_{n+1}$ and with $r+2$ branch
  points.
\end{prop}

\proof We are allowed to replace the field $k$ by an
arbitrary extension $k'$. We may therefore assume that the cover
$f:X\to\PP^1_k$ has generic branch points, see \S \ref{2pt}. (To see
this, let $\HH_{n,r}$ denote the Hurwitz space parameterizing tame
covers of $\PP^1$ of degree $n$ with $r$ branch
points. By \cite{diss}, Proposition 4.2.1, $\HH_{n,r}$ is a smooth
scheme of relative dimension $r$ over $\ZZ$. The cover $f$ of the
proposition corresponds to a $k$-rational point on $\HH_{n,r}\otimes
k_0$. We may thus replace $f$ by the generic cover corresponding to
the connected component of $\HH_{n,r}\otimes k_0\alg$ on which this
point lies.)

Let $z:\Spec k_1\to\PP^1_k$ be a generic geometric point. Let
$h_1:Y_1\to\PP^1_{k_1}$ denote the Galois closure of $f\otimes
k_1$. We identify the Galois group of $h_1$ with $G_1:=A_n$, which we
regard as a subgroup of $G:=A_{n+1}$, in the obvious way. We choose a
point $y_1\in Y_1$ which lies above $z\in\PP^1_{k_1}$. If $p\not=3$
(resp.\ if $p=3$), we let $h_2:Y_2\to\PP^1_{k_1}$ be a cyclic cover of
order $3$ (resp.\ of order $2$), ramified at two points. We identify
the Galois group of $h_2$ with the subgroup $G_2\subset G$ generated
by the $3$-cycle $(n-1,n,n+1)$ (resp.\ by the double transposition
$(1,2)(n,n+1)$). We also choose a point $y_2\in Y_2$ where $h_2$ is
unramified. Let $k':=k_1((t))\alg$. We lift the datum
$(G,G_i,h_i,y_i)$ to a tame $G$-Galois cover $h':Y'\to\PP^1_{k'}$, see
\S \ref{deform}. Let $H\subset G$ be the stabilizer of $1$. Set
$X':=Y'/H$ and let $f':X'\to\PP^1_{k'}$ denote the natural map. By
construction, $f'$ is a tame cover of degree $n+1$, with monodromy
group $A_{n+1}$ and with $r+2$ branch points. The Riemann--Hurwitz
formula shows that $X'$ has genus $g+1$. It remains to prove that $X'$
is generic. 

Let $\X$ be the semistable model of $X'$ over $k_1[[t]]$ originating
from the construction of \S \ref{deform}. If $p\not=3$, then the
special fiber $X_s$ of $\X$ contains the curves $X\otimes k_1$ and
$Y_2$ as irreducible components. These components meet in two points
$x_1,x_2\in X$ with $f(x_1)=f(x_2)=z$. The other components of $X_s$
are copies of $\PP^1_{k_1}$.  If $p=3$, the situation is similar. The
only difference is that $X_s$ contains two copies of $Y_2$ as
irreducible components, corresponding to the orbits of length two of
$(1,2)(n,n+1)$ acting on $\{1,2,\ldots,N\}$.

If $g=0$, the above description of $X_s$ shows that $X'$ has bad
reduction and is therefore a generic curve of genus $1$. Hence we may
assume that $g\geq 1$.  Let $\X\stab$ be the stable model of $X'$ over
$k_1[[t]]$ (i.e.\ the minimal semistable model). The curve $\X\stab$
is obtained from the curve $\X$ by blowing down all components of
$X_s$ of genus $0$ which contain less than three singular points. It
follows that the special fiber of $\X\stab$ is isomorphic to the curve
$X\otimes k_1$, with the points $x_1$ and $x_2$ identified. Let
$\varphi:\Spec k_1[[t]]\to\bar{\M}_{g+1}$ denote the classifying map
for $\X\stab$, and write $\varphi_\eta$ (resp.\ $\varphi_s$) for the
restriction of $\varphi$ to the generic (resp.\ the closed) point of
$\Spec k_1[[t]]$. We have to show that the image of $\varphi_\eta$ is
Zariski dense in $\M_{g+1}$. By our description of the special fiber,
the image of $\varphi_s$ lies on a boundary component $T$ of
$\bar{\M}_{g+1}$ which is isomorphic to $\M_{g,2}$.  Moreover,
$\varphi_s$ corresponds, under this isomorphism, to the $2$-pointed
curve $(X\otimes k_1;x_1,x_2)$ (\cite{KnudsenII}). It follows from
Lemma \ref{2ptlem} that the image of $\varphi_s$ is Zariski dense in
the boundary component $T$ (note that $A_n$ is doubly transitive for
$n\geq 4$ and that the case $n=3$ does not occur for $g\geq 1$). But
$T$ has codimension one in $\bar{\M}_{g+1}$ (\cite{DeligneMumford}).
Therefore, the image of $\varphi_\eta$ is Zariski dense in
$\M_{g+1}$. This concludes the proof of Proposition \ref{indprop2}.
\Endproof

\section{The main result}

\subsection{The case $g=0$} \label{g=0}

Proposition \ref{g=0prop}  implies Part (a) of Theorem \ref{thm0}.  In
the proof we use the following easy lemma.

\begin{lem}\label{ellipticlem}
Let $p\neq 2$. Choose $\lambda \in k-\{0,1\}$ such that
\[
    \Phi(\lambda)=\sum_{i=0}^{(p-1)/2}\binom{(p-1)/2}{i}^2\lambda^i
\]
is nonzero. Then there exists a cover $f:\PP^1_k\to\PP^1_k$ of degree
$p$ branched at $0,1,\lambda,\infty$ with monodromy group $D_p$. For
$x=0,1,\lambda,\infty$, the inverse image of $x$ consists of $(p-1)/2$
ramification points of order two and one unramified point.
\end{lem}

\proof Choose $\lambda \in k-\{0,1\}$ such that $\Phi(\lambda)$ is
 nonzero. This implies that the elliptic curve $E$ defined by
 $y^2=x(x-1)(x-\lambda)$ is ordinary, see e.g.\ \cite{Hartshorne},
 Corollary IV.4.22. Let $E'\to E$ be an \'etale
 isogeny of degree $p$. Taking the quotient on both sides under the
 elliptic involution, we obtain a rational map $f:\PP^1_k\to\PP^1_k$
 of degree $p$ as in the statement of the lemma.  \Endproof

\begin{prop} \label{g=0prop}
  Let $n\geq 3$ and $p\geq 0$. Suppose that
  $n\not=4$ if $p=2$ and $n\geq 7$ if $p=3$. Then there exists
  a tame cover $f:\PP^1_k\to\PP^1_k$ of degree $n$ with monodromy
  group $A_n$, defined over an algebraically closed field $k$ of
  characteristic $p$. More precisely:
\begin{itemize}
\item[(a)] If $p\not=2,3$, we may choose $f$ with $r=n-1$ branch
  points and ramification type $(3,\ldots,3)$.
\item[(b)] Suppose that one of the following holds: 
  \begin{enumerate}
    \item  $p\geq 7$ and $n\geq 6$, 
    \item  $p=3$ and $n\geq 7$, or 
    \item  $p=5$ and $n\geq 7$ is odd. 
  \end{enumerate}
  Then we may choose $f$ with $r=n-1$ branch points and
  ramification type $(\text{$2$-$2$},\ldots,\text{$2$-$2$})$.
\item[(c)] If $p=2$ and $n$ is odd, we may choose $f$ with $r=n-1$
  branch points and ramification type $(3,\ldots,3)$. 
\item[(d)] If $p=2$ and
  $n\geq 6$ is even, we may choose $f$ with $r=n-3$ branch points, and
  ramification type $(5,\text{$3$-$3$},3,\ldots,3)$.
\end{itemize}
\end{prop}  

\proof Suppose first that $p\not=2,3$ and consider the case of triple
ramification. The case $n=3$ is trivial, as we may take for $f$ the
function $f(x)=x^3$. For $n=4$, we set
\begin{equation} \label{333eq}
   f(x) \;:=\; \frac{x^3\,(x-2)}{1-2x} \;=\;
      \frac{(x-1)^3\,(x+1)}{1-2x}+1.
\end{equation}
It is obvious that $f$ has ramification locus $\{0,1,\infty\}$ and
ramification index $3$ in each of these points. (Note that we use
$p\neq 2,3$.) The monodromy group of $f$ is a transitive permutation
group on $4$ letters generated by $3$-cycles, hence isomorphic to
$A_4$. This settles the case $n=4$. The cases $n>4$ follow from
$n=3,4$ by induction, using Proposition \ref{indprop1}. This proves
(a)

We now prove (b). If $p\geq 7$ and $n=6$, we
may use Riemann's Existence Theorem, since in this case $p$ does not
divide the order of $A_6$. For example, one checks using GAP that the
tuple 
\[
((1,2)(3,4); (3,6)(4,5); (2,6)(3,4); (2,6)(3,5); (1,2)(4,6))
\]
generates $A_6$.

Next we handle the case $p\not=2$ and $n=7$. Let
$f_1:\PP^1_k\to\PP^1_k$ be a separable map of degree four branched at
four points $0,1,\lambda, \infty$ with ramification type $(2,2,2,2)$
and monodromy group $G_1$ the dihedral group of order $8$. For $p\neq
2$, the existence of such a cover follows from Riemann's Existence
Theorem. Write $h_1:Y_1\to \PP^1_k$ for the Galois closure of $f_1$
and choose a point $y_1\in Y_1$ above $x=0\in \PP^1_k$. Choose an
embedding of $G_1$ into $A_7$ such that $G_1$ acting on
$\{1,2,\ldots,7\}$ has orbits $\{1,2\}$, $\{3,4,6,7\}$, $\{5\}$. Write
$(g_1,g_2,g_3,g_4)$ for the corresponding tuple of double
transpositions.  Choose $\mu\in k-\{0,1\}$. If $p=5$ we suppose
moreover that $\mu$ is not a primitive $3$rd root of unity. Let
$f_2:\PP^1_k\to\PP^1_k$ be a separable map of degree five branched at
$0,1,\mu,\infty$ with monodromy group the dihedral group of order $10$ and
with ramification type $(\zz,\zz,\zz,\zz)$. If $p\neq 5$ the existence
of $f_2$ follows from Riemann's Existence Theorem. If $p=5$ it follows
from Lemma \ref{ellipticlem}. Write $h_2:Y_2\to\PP^1_k$ for the Galois
closure of $f_2$ and choose a point $y_2\in Y_2$ above $x=0\in
\PP^1_k$. Let $G_2$ denote the Galois group of $h_2$. We embed $G_2$
into $A_7$ such that the action of $G_2$ on $\{1,2,\ldots,7\}$ has
orbits $\{1,2,3,4,5\}$, $\{6\}$, $\{7\}$.  Write
$(g'_1,g'_2,g'_3,g'_4)$ for the tuple of transpositions corresponding
to the cover $f_2$. We may arrange things such that $g_1=g'_1$. One
checks that $A_7$ is generated by $G_1$ and $G_2$. Therefore the
results of \S \ref{deform} imply that there exists a map
$\PP^1_k\to\PP^1_k$ of degree $n=7$ with monodromy group $A_7$ and
ramification $(\zz,\zz,\zz,\zz,\zz,\zz)$. This settles the case
$p\not=2$ and $n=7$.

Part (i) and (iii) of (b) now follow from the cases
$n=6,7$ by induction, using Proposition \ref{indprop1}. In order
to prove Part (ii) of (b) we need an extra construction to handle the case
$n=8$, because the above construction for $n=6$ worked only for $p\geq
7$.

We define (in characteristic $3$) a tame cover $f_1:\PP^1_k\to\PP^1_k$ of
degree $5$ with ramification type $(\zz,\zz,5)$ and Galois group
$D_5$. Since $3$ does not divide the order of the Galois group, we may
use Riemann's Existence Theorem. For example, we may choose $f_1$ such
that it corresponds to the triple $( (1,2)(3,4); (1,5)(2,3);
(1,5,2,4,3) )$.

We also define (in characteristic $3$) a cover 
\[
  f_2(x)\;=\; \frac{x^5(x+1)}{x-1}\;=\; \frac{(x^2-x-1)^2(x^2+1)}{x-1}+1.
\]
The cover $f_2:Y_2\to\PP^1_k$ is branched at $x=0,1,\infty$, with
ramification type $(5,\zz,5)$. One easily checks that this is the only
such cover in characteristic $3$, up to isomorphism. Let $G_2\subset
A_6$ be its monodromy group. One checks, for example using GAP, that
$G_2$ is either isomorphic to $\PSL_2(5)\cong A_5$ or to $A_6$. 

We claim that $G_2\simeq\PSL_2(5)$. (One does not actually need this
in what follows.) It follows from \cite[Proposition
7.4.3]{SerreTopics} that there is a unique $\PSL_2(5)$-Galois cover
$g:Y\to\PP^1_{\bar{\QQ}}$ branched at three points with ramification
type $(5,\zz,5)$ over $\bar{\QQ}$, up to isomorphism. In fact, this
cover may be defined over $\QQ$. It follows from the discussion in
\cite[\S 2.4]{RRR} that $g$ has good reduction to characteristic
$p=3$. Here one uses that $p$ exactly divides the order of
$\PSL_2(5)$. Since $f_2$ is the unique cover (up to isomorphism) of
degree $6$ with ramification type $(5,\zz,5)$ in characteristic $3$,
we conclude that the reduction of $g$ to characteristic $3$ is
isomorphic to $f_2$. This shows that $G_2\cong\PSL_2(5)$. 

Let $g_1\in G_2$ (resp.\ $g_3\in G_2$) be canonical generators of
inertia at the point $x=0$ (resp.\ $x=\infty$), with respect to some
fixed $5$th root of unity. Let $f_3=f_1$. For $i=1,3$, we choose an
embedding of the monodromy group $G_i$ of $f_i$ into $\PSL_2(5)$ by
identifying $G_i$ with the normalizer in $\PSL_2(5)$ of the subgroup
generated by $g_i$. We may choose this identification such that the
canonical generator of inertia of the ramification point of order $5$
of $f_i$ is $g_i^{-1}$. Using the results of \S \ref{deform}, we may
patch the covers $f_1, f_2, f_3$, yielding a cover
$f:\PP^1_k\to\PP^1_k$ of degree $6$ with monodromy group $G_2\subset S_6$
and ramification type $(\zz,\zz,\zz,\zz,\zz)$.

We now apply the construction of Proposition \ref{indprop1} to
$f$. Note that it does not matter for the construction that the
monodromy group of $f$ is $\PSL_2(5)$ and not $A_6$. This yields a
cover $f':\PP^1_k\to\PP^1_k$ of degree $8$ branched at $7$ points of
type $(\zz,\ldots,\zz)$. Note that the monodromy group $G'$ of $f'$ is
a transitive group on $8$ letters which contains $\PSL_2(5)$ and is
contained in $A_8$. One checks, e.g.\ using GAP, that the only such
group is $A_8$ itself. This proves (b) for $p=3$ and $n=8$. Now Part
(ii) of (b) follows from the cases $n=7,8$ by induction, using
Proposition \ref{indprop1}. This completes the proof of (b).

Let us now prove (c) and (d). We suppose that $p=2$. The case $n=3$ is again
trivial; the general case of (c) follows by induction, using Proposition
\ref{indprop1}. To prove (d), we have to start the induction with
$n=6$. We define (in characteristic $2$)
\begin{equation} \label{533eq}
  f(x) \;=\; \frac{x^3\,(x^3+x^2+1)}{x+1} 
       \;=\; \frac{(x^2+x+1)^3}{x+1}+1.
\end{equation}
One checks that $f$ is a tame cover of degree $6$, branched at
$\{0,1,\infty\}$ and with ramification type $(3,\text{$3$-$3$},5)$. A
standard result in group theory (see e.g.\ \cite{Huppert}, p.\ 171)
shows that the monodromy group of $f$ is isomorphic to $A_6$. This 
settles the case $(p,n)=(2,6)$. Part (d) of the proposition follows
now from the cases $n=6$ by induction, using Proposition \ref{indprop1}.
This completes the proof of Proposition \ref{g=0prop}.
\Endproof


The statements of Proposition \ref{g=0prop} is probably not
optimal. In the next proposition we list all the cases where we can
positively exclude the existence of a tame rational function on
$\PP^1$ with alternating monodromy. 

\begin{prop} \label{g=0prop2}
\begin{enumerate}
\item
  Let $f:\PP^1_k\to\PP^1_k$ be a tame rational function of degree
  $n=3$ or $n=4$, with monodromy group $A_n$. Then the ramification
  type of $f$ is $(3,\ldots,3)$ (with $n-1$ branch points).
\item
  Suppose that $p=2$. Then there does not exist a tame rational
  function $f:\PP^1_k\to\PP^1_k$ of degree $4$ and with monodromy
  group $A_4$.
\item
  Suppose that $p=3$. Then there does
  not exist a tame rational function $f:\PP^1_k\to\PP^1_k$ of degree
  $n\leq 5$ with monodromy group $A_n$. 
\end{enumerate}
\end{prop}

\proof The proof of (i) is trivial.  To prove (ii), suppose that there
does exist a tame cover $f:\PP^1_k\to\PP^1_k$ of degree $4$ with
monodromy group $A_4$. By (i), such a cover would be of type
$(3,3,3)$. Such a cover would also lift to characteristic $0$, in a
unique way. It is easy to check, using the rigidity criterion
(\cite[Section 7.3]{SerreTopics}), that there is a unique $A_4$-cover
of $\PP^1$ with ramification type $(3,3,3)$ in characteristic
zero. Therefore, the cover we are looking for would be the reduction
of the cover \eqref{333eq} to characteristic $2$. However, it is easy
to check that the cover \eqref{333eq} has bad reduction to
characteristic $2$. This gives the desired contradiction and proves
(ii).

To prove (iii), suppose that $p=3$. The cases $n=3,4$ are already
excluded by (i). It is easy to see, using the Riemann--Hurwitz
formula, that a tame cover $f:\PP^1_k\to\PP^1_k$ of degree $5$ has
ramification type $(5,5)$, $(5,\zz,\zz)$ or $(\zz,\zz,\zz,\zz)$. In
the first case the monodromy group is cyclic of order $5$, in the
second case a dihedral group of order $10$.  So it remains to rule out
the third case. Let $(g_1,g_2,g_3,g_4)$ be a $4$-tuple of double
transpositions in $A_5$ with $g_1g_2g_3g_4=1$. We claim that the
elements $g_1,\ldots,g_4$ do not genererate $A_5$. To prove the claim,
let $h:=g_1g_2$ and $G$ the group generated by $g_1,\dots,g_4$. If
$h=1$ then $G$ is generated by two elements of order $2$; it is then
isomorphic to $\ZZ/2$, $S_3$ or $D_5$, but not to $A_5$. If $h$ is a
$3$-cycle, then the subgroups $G_1:=\gen{g_1,g_2}\subset A_5$ and
$G_2:=\gen{g_3,g_4}$ are both isomorphic to $S_3$ and contain the same
subgroup of order $3$. It follows easily that $G_1=G_2$ and hence that
$G\cong S_3$. Finally, if $h$ is a $5$-cycle, then $G_1$ and $G_2$ are
dihedral groups of order $10$ and contain the same cyclic group of
order $5$. As before, it follows that $G=G_1=G_2\cong D_5$. This
proves the claim, and completes the proof of Proposition
\ref{g=0prop2}.  \Endproof

\begin{rem} \label{g=0rem}
 The following cases are left open by Proposition \ref{g=0prop} and
 Proposition \ref{g=0prop2}:
 \begin{enumerate}
 \item
   $p=2$, $n\geq 6$ even, and ramification type $(3,\ldots,3)$,
 \item
   $p=3$, $n=6$,
 \item
   $p=5$, $n\geq 6$ even, and ramification type $(\zz,\ldots,\zz)$.
\end{enumerate}
\end{rem}

\subsection{}

The following theorem is a more precise version of  Theorem \ref{thm0}
of the introduction. 

\begin{thm}\label{thm1} 
  Let $g\geq 0$ and $n\geq 3$. Let $X$ be a generic curve of genus
  $g$, defined over an algebraically closed field $k$ of
  characteristic $p\geq 0$. Then there exists a tame rational function
  $f:X\to\PP^1_k$ of degree $n$ with monodromy group $A_n$ in each of
  the following cases.
\begin{itemize}
\item[(a)] If $p\neq 2,3$ and $n\geq\max\,(g+3,2g+1)$, we may choose $f$ with
  $r=g+n-1$ branch points and ramification type $(3,\ldots,3)$.
\item[(b)] Suppose that one of the following holds: 
  \begin{enumerate}
    \item  $p\geq 7$ and $n\geq \max\,(6+g,2g+1)$, 
    \item  $p=3$ and $n\geq \max\,(7,6+g,2g+1)$, or 
    \item  $p=5$, $n\geq \max\,(7+g,2g+1)$ and $n+g$ is odd. 
  \end{enumerate}
  Then we may choose $f$ with $r=g+n-1$ branch points and
  ramification type $(\text{$2$-$2$},\ldots,\text{$2$-$2$})$.
\item[(c)] If $p=2$, $n\geq\max\,(g+3,2g+1)$ and $n+g$ is odd, then we
  may choose $f$ with $r=g+n-1$ branch points and ramification type
  $(3,\ldots,3)$.
\item[(d)] If $p=2$, $n\geq\max\,(g+6,2g+3)$ and
  $n+g$ is even, we may choose $f$ with $r=g+n-3$ branch points, and
  ramification type $(5,\text{$3$-$3$},3,\ldots,3)$.
\end{itemize}
\end{thm}

\proof We start with two observations. First, to prove the theorem we
may extend the base field $k$ by an arbitrary algebraically closed
field extension $k'/k$. Indeed, if we can show that the curve
$X\otimes k'$ admits a rational function with certain properties, then
a standard specialization argument yields the existence of a rational
function on $X$ with the same properties. Second, by the
irreducibility of $\M_g$, a generic curve of a given genus is unique,
up to isomorphism and extension of the base field. We may hence speak
about {\em the} generic curve of genus $g$. These two observations
will allow us to prove the theorem by induction on the pair $(g,n)$,
using Proposition \ref{indprop2}.

We will discuss the induction procedure in detail for Part
(a). For Part (b)--(d), we will only give the necessary
modifications. 

Suppose that $p\not=2,3$. For $g=0$, the statement of Theorem
\ref{thm1} (a) is equal to the statement of Proposition \ref{g=0prop}
(a). Fix an integer $g'\geq 1$ and suppose that we have already proved
Theorem \ref{thm1} (a) for all pairs $(g,n)$ with $g<g'$. We have to
show that Theorem \ref{thm1} (a) holds for all pairs $(g',n')$ with
$n'\geq\max\,(g'+3,2g'+1)$. Write $g'=g+1$ and $n'=n+1$. Let $X$ be
the generic curve of genus $g$ over $k$. By the induction hypothesis,
there exists a tame rational function $f:X\to\PP^1$ of degree $n$ with
alternating monodromy, $r=g+n-1$ branch points and ramification type
$(3,\ldots,3)$. We also have $n>2g+1$ and hence $r>3g$. Therefore, by
Proposition \ref{indprop2}, the generic curve $X'$ of genus $g'$
admits a tame rational function $f':X'\to\PP^1$ of degree $n'$, with
monodromy group $A_{n'}$ and ramification type $(3,\ldots,3)$. In
other words, Theorem \ref{thm1} (a) holds for the pair $(g',n')$. This
completes the proof of (a).

The proof of (b) is almost the same. Note that for $g=0$ the statement
of Theorem \ref{thm1} (b) reduces again to the statement of
Proposition \ref{g=0prop} (b). The only problem occurs for $p=3$ and
the pair $(g',n')=(1,7)$, because Proposition \ref{g=0prop} (b) says
nothing about the case $(g,n)=(0,6)$ (see also Proposition
\ref{g=0rem} (ii)). However, in the proof of Proposition \ref{g=0prop}
(b) we did construct a tame rational function $f:\PP^1\to\PP^1$ of
degree $6$, with $5$ branch points and ramification type
$(\zz,\ldots,\zz)$. We showed that its monodromy group is isomorphic
to $\PSL_2(5)\cong A_5$. If we apply the construction underlying the
proof of Proposition \ref{indprop2} to this cover, then we obtain a
tame rational function $f':X'\to\PP^1$ on the generic curve of genus
one of degree $7$, with ramification type $(\zz,\ldots,\zz)$. Its
monodromy group is a transitive subgroup of $A_7$ and contains a
subgroup isomorphic to $A_5$. One checks, e.g.\ using GAP, that this
group must be $A_7$. This settles the case $(p,g,n)=(3,1,7)$
and completes the proof of (b).  

The proof of (c) and (d) is again similar. The only difference is that
we have a different bound on $n$ in (d). However, the induction
procedure used in the proof of (a) and (b) goes through. 
\Endproof

\subsection{Open problems in characteristic $2$ and $3$} \label{OpenProblems}

As we already mentioned in the introduction, our main result (Theorem
\ref{thm1}) is optimal for $p\not=2,3$ (in a certain sense). However,
in characteristic $2$ and $3$ there remain a number of cases which we
could not handle with the methods of the present paper. We leave
these remaining cases as problems for the interested reader.

\begin{prob} \label{problem1}
  Let $X$ be the generic curve of genus $1$ in characteristic
  $2$. Show that there exists a tame rational function $f:X\to\PP^1$
  of degree $5$ with monodromy group $A_5$  and ramification type
  $(3,3,3,3,3)$. 
\end{prob}

Let $g:E\to\PP^1$ be the tame cyclic Galois cover of degree $3$ with
three branch points. Note that, in characteristic $2$, the curve $E$
is the unique supersingular elliptic curve. (In characteristic $0$,
$E$ is an elliptic curve with complex multiplication by
$\ZZ[\zeta_3]$.) Applying the construction underlying the proof of
Proposition \ref{indprop1}, we obtain a tame cover $f:X\to\PP^1$ of
degree $5$ with monodromy group $A_5$, ramification type $(3,3,3,3,3)$
and with generic branch points. By Riemann-Hurwitz, $X$ has genus
$1$. We believe but have not been able to show that $X$ is
generic. The problem is the following. By construction, the cover
$f:X\to\PP^1$ is the generic fibre of a family of covers over a
two-dimensional base. On a sublocus of codimension one of the base,
this family degenerates to the cover $g:E\to\PP^1$. In particular,
over this sublocus the top curve of the cover is isotrivial. But it is
not clear whether $X$ is isotrivial or generic.

Note that the characteristic $0$ version of Problem \ref{problem1} is
solved in the paper \cite{FriedKlaKop00}, using the same construction
as above. The proof uses the braid action to show the $X$ is not
isotrivial. In fact, the case $(g,n)=(1,5)$ is the only case where the
braid action enters into the proof of the main result of
\cite{VoelkleinMagaard03} in an essential way. 

\begin{rem}
  At the moment, we do not know whether the generic curve of genus $g$
  in characteristic $2$ admits a tame rational function of degree $n$
  with alternating monodromy, in each of the following cases:
  \begin{enumerate}
  \item
    $(g,n)=(1,5)$,
  \item
    $g\geq 3$ is odd and $n=2g+1$, 
  \item
    $g\geq 2$ is even and $n=2g+2$.
  \end{enumerate}
  If Problem \ref{problem1} (which corresponds to Case (i)) had a
  positive solution, then also Case (ii) and (iii) would be solved, as
  one can see from the induction procedure in the proof of Theorem
  \ref{thm1}. If this were the case, then we could omit the condition
  `$n+g$ is odd' in the statement of Theorem \ref{thm1} (c) (but we
  would have to add the condition `$(g,n)\not=(0,4)$'). This would
  make Part (d) of Theorem \ref{thm1} superfluous, and would give an
  `optimal' result even in characteristic $2$.

  In characteristic $3$, the existence of tame rational functions with
  alternating monodromy is open in the following (finite) list of
  cases:
  \[
     (g,n) \;=\; (1,5),(1,6),(2,5),(2,6),(2,7),(3,7),(3,8),(4,9).
  \]
\end{rem}

\end{document}